\documentclass{amsart}
\usepackage{amssymb}
\usepackage{amscd}
\usepackage{amssymb}
\usepackage{amscd}
\usepackage{epsfig}
\title[Double Shuffle Relations]{Regularization and generalized double shuffle relations for $p$-adic multiple zeta values}
\author{Hidekazu Furusho}
\author{Amir Jafari}
\address{Graduate School of Mathematics, Nagoya University,
Nagoya, 464-8602, Japan}
\email{furusho@math.nagoya-u.ac.jp}
\address{Math Departments, Duke University, Durham, NC, 27708, USA}
\email{amir@math.duke.edu}
\newcommand{\Q}{\bold Q}
\newcommand{\Z}{\bold Z}
\newcommand{\C}{\bold C}
\newcommand{\N}{\bold N}

\newcommand{\To}{\longrightarrow}
\newcommand{\M}{{\mathcal M}_{0,N+3}}

\newcommand{\Mi}{{\mathcal M}_{0,4}}
\newcommand{\MM}{\overline{\mathcal M}_{0,N+3}}

\newtheorem{thm}{Theorem}[section]
\newtheorem{lem}[thm]{Lemma}
\newtheorem{cor}[thm]{Corollary}
\newtheorem{prop}[thm]{Proposition}

{\theoremstyle{definition} }
{\theoremstyle{definition} \newtheorem{defn}[thm]{Definition}}

{\theoremstyle{remark} }

\numberwithin{equation}{section}

\begin{document}

\bibliographystyle{amsalpha+}
\maketitle

\begin{abstract}
  We will introduce a regularization for $p$-adic multiple zeta
  values and show that the generalized double shuffle
  relations hold. This settles a question raised by Deligne, given as a project
  in Arizona winter school 2002. Our approach is to use the theory of Coleman
  functions on the moduli space of genus zero curves with marked
  points and its compactification. The main ingredients are the
  analytic continuation of Coleman functions to the normal bundle
  of divisors at infinity and definition of a special
  tangential base point on the moduli space.

\end{abstract}
\tableofcontents
\setcounter{section}{-1}
\section{Introduction}

 This paper is a continuation of \cite{BF}. Let $p$ be a prime
 number. The first author in \cite{Fu1} introduced $p$-adic MZV's \footnote{MZV stands for multiple zeta
 value.} for admissible indices $(n_1,\dots,n_m)\in\Z^m_{>0}$ that is when $n_m>1$.
 In \cite{BF} the double shuffle relations for these $p$-adic
 MZV's were proved.
 In this work we extend their result.
 By using the moduli space of genus zero curves with marked
 points and its stable compactification we introduce a series regularization of $p$-adic MZV's for
 non-admissible indices and then prove a generalized double shuffle relation.
 This relation is an extension of the double shuffle relation in
 \cite{BF}. This includes a comparison between the integral
 regularization and the series regularization of $p$-adic MZV's.

Let us first review the complex story. Recall that for positive
 integers $n_1,\dots,n_m$ the (complex) MZV is defined by:
\begin{equation}\label{MZV}
\zeta(n_1,\dots,n_m)=\sum_{0<k_1<\dots<k_m}\frac{1}{k_1^{n_1}\dots
k_m^{n_m}}.
\end{equation}
They were studied first by Euler for $m=1$ and $m=2$. It is easy
to see that the series is convergent if and only if $n_m>1$. The
{\bf double shuffle relations} consists of {\bf series shuffle
relations} and {\bf integral shuffle relations}. Both of them are
product formulae between MZV's. The simple example of the series
shuffle relation is
$$\zeta(n_1)\zeta(n_2)=\zeta(n_1,n_2)+\zeta(n_1+n_2)+\zeta(n_1+n_2).$$
It is easily obtained from the expression  \eqref{MZV} and can be
generalized in a similar way to other MZV's. The simple example of
the integral shuffle relation is
$$\zeta(n_1)\zeta(n_2)=\sum_{i=0}^{n_1-1}\binom{n_2-1+i}{i}\zeta(n_1-i,n_2+i)+\sum_{j=0}^{n_2-1}\binom{n_1-j+1}{j}\zeta(n_2-j,n_1+j).$$
This follows from the iterated integral expression for MZV's.
Using these formulae we can get many relations among the MZV's
however the double shuffle relations are not enough to capture all
the relations between MZV's. There are two regularization of the
MZV's for non-admissible indices: the {\bf series regularization},
which extends the validity of the series shuffle relation and the
{\bf integral regularization} which extends the validity of the
integral shuffle relation. The several variable MPL\footnote{MPL
stands for multiple polylogarithm.}:
\begin{equation}\label{MPL}
Li_{n_1,\dots,n_m}(x_1,\dots,x_m)=\sum_{0<k_1<\dots<k_m}\frac{x_1^{k_1}\cdots
x_m^{k_m}}{k_1^{n_1}\dots k_m^{n_m}}
\end{equation}
with complex variables $z_1,\dots,z_m$ is the device to construct
the series shuffle regularization. It is clear that when the
indices are admissible the limit of \eqref{MPL} when $z_i$'s
approach 1 is $\zeta(n_1,\dots,n_m)$. In general one can show that
$Li_{n_1,\dots,n_m}(1-\epsilon,\dots,1-\epsilon)=\sum_{i=0}^N
a_i(\epsilon)\log^i \epsilon$, where $a(\epsilon)\in
\C[[\epsilon]]$. The {\bf series regularized value}
$\zeta^S(n_1,\dots,n_m)\in\C[T]$ is by definition the polynomial
$\sum_{i=0}^N a_i(0)T^i$ where $T=\log \epsilon$.
\\
The one variable MPL:
\begin{equation}\label{1MPL}
Li_{n_1,\dots,n_m}(z)=Li_{n_1,\dots,n_m}(1,\dots,1,z)=\sum_{0<k_1<\dots<k_m}\frac{z^{k_m}}{k_1^{n_1}\dots
k_m^{n_m}}
\end{equation}
with one complex variable $z$ is the device to construct the
integral regularization. It is clear that when the indices are
admissible the limit of \eqref{1MPL} when $z$ approaches 1 is
$\zeta(n_1,\dots,n_m)$. It can be shown that:
$Li_{n_1,\dots,n_m}(1-\epsilon)=\sum_{i=0}^M
b_i(\epsilon)\log^i\epsilon$, with
$b_i(\epsilon)\in\C[[\epsilon]]$. The {\bf integral regularized
value} $\zeta^I(n_1,\dots,n_m)\in\C[T]$ is by definition the
polynomial $\sum_{i=0}^M b_i(0)T^i$ where $T=\log \epsilon$. There
is a comparison relation between these two regularization that we
now describe. Let ${\mathbb L}$ be the $\C$ linear map from the
polynomials $\C[T]$ to itself defined via the generating function:
\begin{equation}\label{gen}
\sum_{n=0}^{\infty}{\mathbb
L}(T^n)\frac{u^n}{n!}=\mbox{exp}\left(-\sum_{n=1}^{\infty}\frac{\zeta^I(n)}{n}u^n\right).\footnote{There
is a sign error in \cite{G1} for this formula.}
\end{equation}
The {\bf regularization relation} \cite{IKZ} asserts that
\begin{equation}\label{regg}
\zeta^S(n_1,\dots,n_m)={\mathbb L}(\zeta^I(n_1,\dots,n_m)).
\end{equation}
The {\bf generalized double shuffle relation} is in fact three
types of relations: the series shuffle relation for series
regularized MZV, the integral relation for integral shuffle
regularized MZV and the relation \eqref{regg} that gives the
comparison between these two regularization. This is conjectured
to capture all the possible relations between the MZV's (cf.
\cite{R}).

Now we explain the $p$-adic story. The one variable MPL in
\eqref{1MPL} has a $p$-adic analogue as a Coleman function on
${\mathbb P}^1(\C_p)-\{0,1,\infty\}$. This function depends on a
choice of the branch of $p$-adic logarithm. We denote this
function by $Li_{n_1,\dots,n_m}^a(z)$, where $a\in \Q_p$ is the
value of the $p$-adic log at $p$. In \cite{Fu1} the $p$-adic MZV
$\zeta_p(n_1,\dots,n_m)\in\Q_p$ for an admissible index, is
defined as a certain limit of this Coleman function as $z$
approaches 1. It is shown in loc. cit. that this limit exists and
is independent of the choice of the branch of $p$-adic logarithm
$a$ made before. Using the language of the tangential base point,
it can be explained as follows. The function
$Li_{n_1,\dots,n_m}^a(z)$ is a Coleman function on ${\mathbb
P}^1-\{0,1,\infty\}$. First analytically continue it to the
punctured tangent plane at the point $z=1$. Let $z$ denote the
canonical parameter of ${\mathbb P}^1-\{0,1,\infty\}$. If we let
$t=1-z$ denote the local parameter of this tangent plane we get an
element of $\Q_p[T]$ for $T=\log^at$.  We want to stress that this
polynomial is independent of the choice of $a$ and we define it by
the {\bf integral regularized $p$-adic MZV}
${\zeta_p}^I(n_1,\dots,n_m)$. When $n_m>1$ this is a constant
polynomial, cf. \cite{Fu1}.

 We will use a method inspired by \cite{G1} to
define series regularized $p$-adic MZV. We use a $p$-adic analogue
of the several variable MPL in \eqref{MPL}.
 It is denoted by $Li_{n_1,\dots,n_m}^a(x_1,\dots,x_m)$ where $a$ is the value of the branch of $p$-adic logarithm at $p$.
It is a Coleman function on ${\mathcal M}_{0,m+3}$, the moduli
space of curves of type $(0,m+3)$. We identify this space with
${\mathbb A}^m-\{x_1\cdots x_m=0\quad \mbox{or}\quad 1-x_i\dots
x_j=0\quad i\leqslant j\}.$ The line $(1-t,1-t,\dots,1-t)$ in ${\mathcal
M}_{0,m+3}$ when $t$ approaches 0 intersects a unique divisor
$D_0$ of the stable compactification $\overline{{\mathcal
M}}_{0,m+3}$ at a point $R$. Let $L$ be the punctured line above
$R$ in the normal bundle of $D_0$ minus the zero section. The
{\bf series regularized $p$-adic MZV} is defined by:
$${\zeta}^{S}_p(n_1,\dots,n_m):=Li_{n_1,\dots,n_m}^{a,(D_0)}(\bar{x}_1,\dots,\bar{x}_m)|_L.$$
where $Li_{n_1,\dots,n_m}^{a,(D_0)}(\bar{x}_1,\dots,\bar{x}_m)$ is
the analytic continuation of MPL to the normal bundle of $D_0$
minus the zero section and restricted to the open part of the
divisor $D_0$ outside the other divisors.
This $\zeta_p^{S}(n_1,\dots,n_m)$ is a polynomial in $\Q_p[T]$ where
$T=\log^a t$. It a priori depends on $a$.

The main result of the paper is the following:
\begin{thm}\label{main}
\begin{enumerate}
\item The series regularized MZV
$\zeta_p^{S}(n_1,\dots,n_m)\in\Q_p[T]$ is well-defined, that is,
it is independent of the choice of the branch parameter $a\in\bold Q_p$.

 \item If $n_m>1$ then ${\zeta}^S_p(n_1,\dots,n_m)$ is
constant and is equal to the $p$-adic MZV $\zeta_p(n_1,\dots,
n_m)$ in \cite{Fu1}.

\item The generalized double shuffle relation holds: i.e.
$\zeta_p^S(n_1,\dots,n_m)$ satisfies the series shuffle relation
and $\zeta^I_p(n_1,\dots,n_m)$ satisfies the integral shuffle
relation. Furthermore the regularization relation:
$$\zeta_p^S(n_1,\dots,n_m)={\mathbb
L}_p(\zeta_p^I(n_1,\dots,n_m))$$ holds, where $\mathbb L_p$ is
defined by analogous generating series in \eqref{gen} where we
replace $\zeta^I(n)$ by $\zeta_p^I(n)$.
\end{enumerate}
\end{thm}

This theorem is used together with some results of Racinet and the
first author to prove that

\begin{thm}
Deligne's $p$-adic MZV's satisfy the generalized double shuffle
relations. \end{thm}

The definition of these $p$-adic MZV's is recalled in section
\ref{deligne}. This definition was suggested in 2002 Arizona
winter school and the theorem above solves the project proposed by
Deligne.

{\bf Acknowledgments.} The authors thank H. Nakamura for giving
some useful references and answering some of their questions. H.F.
learned lots about MZV's from M.Kaneko.
He also thanks the Duke university for its generosity and support.

\section{Review of Coleman functions}

We recall some definitions and properties of Coleman function and
tangential base points as developed in \cite{B1},\cite{B2} and
\cite{BF}. We fix a branch of $p$-adic logarithm $\log$ with value
$a\in\Q_p$ at $p$ for the rest of this paper.

Let $X$ be a smooth variety over $K$, a finite extension of
$\Q_p$. Let $\mathcal{NC}(X)$ denote the category of unipotent
vector bundles on $X$ with flat connections. This forms  an
neutral tannakian category and any point $x\in X(K)$ defines a fiber
functor $\omega_x$ from $\mathcal{NC}(X)$ to the category $Vec_K$
of finite dimensional $K$-vector spaces (cf.\cite{D}). In \cite{Vo}, Vologodsky
has constructed a canonical system (after fixing a branch of
$p$-adic logarithm) of isomorphism
$a_{x,y}^{X}:\omega_x\To\omega_y$ for any pair of points in
$X(K)$. The properties of these isomorphism are summarized in
\cite{B2}\S2. Following \cite{B2}, an abstract Coleman function
is a
triple $(M,y,s)$ where $M\in \mathcal{NC}(X)$ and $y$ is a
collection of $y_x\in M_x$ for all $x\in X(L)$ for any finite
extension $L$ of $K$ and $s\in \mbox{Hom}_{\mathcal
O_X}(M,\mathcal{O}_X)$ is a section. This data must satisfy:
\begin{itemize}
\item For any two points $x_1,x_2\in X(L)=X_L(L)$ we have
$a_{x_1,x_2}^{X_L}(y_{x_1})=y_{x_2}$. \item For any field
homomorphism $\sigma:L\To L'$ that fixes $K$ and $x\in X(L)$ we
have: $\sigma(y_x)=y_{\sigma(x)}$.
\end{itemize}
There is a natural notion of morphism between the abstract Coleman
functions. The connected component of an abstract Coleman function
is called a Coleman function. A Coleman function can be considered
as a function on $X(\overline{K})$ by assigning to $x$ the value
$s(y_x)$. This is indeed a locally analytic function. We will use
both approaches to Coleman functions in this paper. The set of
Coleman functions with values in ${\mathcal O}_X$ is a ring which
we denote by $\mbox{Col}^a(X)$. Here $a\in \Q_p$ is the value of
the chosen branch of the $p$-adic logarithm at $p$.
\\
Let $X$ be a smooth $\mathcal O_K$ scheme and $D=\sum D_i$ be a
divisor with relative normal crossings over $\mathcal O_K$, with
$D_i$'s smooth and irreducible over $\mathcal O_K$. In \cite{BF} a
tangential morphism $Res_{D,I}:\mathcal{NC}((X-D)_K)\To
\mathcal{NC}({\mathcal N}_I^{00})$ was constructed. Here
${\mathcal N}_I^{00}$ is the normal bundle of $D_I=\cap_{i\in I}
D_i$ minus the normal bundles of $D_{I-\{i\}}$ for all $i\in I$
(the normal bundle ${\mathcal N}_\emptyset$ is considered as the
zero section) , and then restricted to $D_I-\cup_{j\not\in I}
(D_j\cap D_I)$. Let $\kappa$ be the residue field of $\mathcal
O_K$. It was shown in \cite{BF} that if the Frobenius endomorphism
of $(X,D)_{\kappa}$ locally lifts to an algebraic endomorphism of
$(X,D)$ then this morphism respect the action of the Frobenius
endomorphism. Indeed by 
\cite{S}, \cite{CLS} the categories $\mathcal{NC}(X-D)$ and
$\mathcal{NC}({\mathcal N}_I^{00})$ are isomorphic to the
categories of the unipotent isocrystals
$\mathcal{NC}^{\dagger}((X-D)_{\kappa})\otimes K$ and
$\mathcal{NC}^{\dagger}(({\mathcal N}_I^{00})_{\kappa})\otimes K$
on the reductions $(X-D)_{\kappa}$ and $({\mathcal
N}_I^{00})_{\kappa}$ and therefore admit a natural action of the
Frobenius endomorphism. Choose a point $\tilde{t}\in ({\mathcal
N}^{00}_I)_{\kappa}(\overline{\kappa})$ which is the reduction of
a point $t\in {\mathcal N}_I^{00}(L)$ for some extension $L$ of
$K$. The point $\tilde{t}$ defines a fiber functor
$\omega_{\tilde{t}}$ from $\mathcal{NC}^{\dagger}((X-D)_{\kappa})$
to $Vec_L$, which is Frobenius invariant if we take a high power
of the Frobenius. Then following \cite{B1} for any point
$\tilde{x}\in (X-D)_{\kappa}(\overline{\kappa})$, which is the
reduction of $x\in {\mathcal N}_I^{00}(L)$, we get a canonical Frobenius
invariant isomorphism
$\tilde{a}_{\tilde{x},\tilde{t}}:\omega_{\tilde{x}}\To
\omega_{\tilde{t}}$. 
The above categorical equivalence gives
an isomorphism $a_{x,t}:\omega_x\To \omega_t$. Now for any $x'\in
(X-D)(L)$ and $t'\in {\mathcal N}_I^{00}(L)$ we define :
$$a_{x',t'}=a_{x',x}\circ a_{x,t}\circ a_{t,t'}.$$
This is independent of the choice of $x$ and $t$. Using this we
can define a way to extend a Coleman function (with values in
$\mathcal O_X$), $(M,\nabla,y,s)$, to a Coleman function on
$\N_I^{00}$ as follows: the connection will be
$Res_{D,I}(M,\nabla)$. The section $s$ will be the induced section
via the functor $Res_{D,I}$. We define $y_t=a_{x,t}(y_x)$ for some
$x\in (X-D)(K)$. We refer the reader to \cite{BF}\S 2 for more
details.

\section{The moduli space ${\mathcal M}_{0,N+3}$ and its compactification}

In this section we give a quick review on some basic properties of
the moduli space $\M$ of genus zero curves with $N+3$ distinct
marked points and its stable compactification $\MM$. The basic
references are \cite{GHP}, \cite{GM} and \cite{Man}.\par
The moduli space $\M$ can be identified with
$${\mathbb A}^N-\left\{x_n=0, \quad \prod_{i=n}^m x_i=1,\quad
1\leqslant n\leqslant m\leqslant
N\right\}.$$ The identification is given by sending
$(x_1,\dots,x_N)$ to the $N+3$ marked points on ${\mathbb P}^1$
given by
$(0, x_1\cdots x_N, x_2\cdots x_N,\dots, x_N, 1, \infty).$
Note that with this identification we have canonical coordinates $x_1,\dots ,x_N$ on $\M$.

We need to work with $\MM$, the stable compactification of this
moduli space. There is a very concrete description of this space
in \cite{GHP} that we now recall. Let $V_N$ be the set of all
distinct ordered 4-tuples of $\{1,\dots,N+3\}$. There is an
embedding:
$$r:\M\hookrightarrow {\mathbb A}^{V_N}$$
given by sending $(P_1,\dots,P_{N+3})$ to all cross ratios of
4-tuple of points. To normalize the cross ratio we recall that
$r(0,\infty,1,x)=x$. Let $\lambda_v$ for $v\in V_N$ be the
coordinates of ${\mathbb A}^{V_N}$.  The image variety will be
given by the following equations:
$$\lambda_{v_1v_2v_3v_4}\lambda_{v_2v_1v_3v_4}=1$$
$$\lambda_{v_1v_2v_3v_4}=1-\lambda_{v_2v_3v_4v_1}$$
$$\lambda_{v_1v_2v_4v_5}\lambda_{v_1v_2v_3v_4}=\lambda_{v_1v_2v_3v_5}$$
for all distinct 5-tuples $v_1,v_2,v_3,v_4,v_5$ in
$\{1,\dots,N+3\}$. Now the compactification $\MM$ is obtained
simply by taking the closure of this variety inside $({\mathbb
P}^1)^{V_N}$. This means that we homogenize the equations by
letting $\lambda_v=\frac{a_v}{b_v}$.

To give a natural stratification of $\MM$ which is of
combinatorial origin we need to recall the notion of stable
labeled trees. A stable tree is a tree such that each of its
vertices have valency $\geqslant 3$. This implies that the tails of this
tree are open in one side. A stable $(N+3)$-labeled tree is a tree
with $N+3$ tails labeled by distinct labels from the set
$\{1,\dots,N+3\}$. The tail of a tree is a one sided open edge.
The stability condition means that the valency of all the vertices
are $\geqslant 3$. Any stable curve of genus zero with $N+3$ marked
points defines a stable $(N+3)$-labelled tree. This construction is
standard and for details consult the above references.

Any vertex $t$ of a stable $(N+3)$-labelled tree $T$ defines an
equivalence relation $\sim_t$ on the set $\{1,\dots,N+3\}$ which
can be identified with the set of tails, as follows:
$i\sim_t j$ if either the corresponding tails have a common vertex
or there is a path from $i$ to $j$ in $T$ that avoids $t$.\par

To each stable $(N+3)$-labeled tree $T$ we associate a closed smooth
subvariety of $\MM$, denoted by $D(T)$. In coordinates it is
defined by:
$$\lambda_v=0\quad \mbox{for all}\quad v\in V(T)$$
where $V(T)$ is the subset of $V_N$ of those quadruples
$(v_1,v_2,v_3,v_4)$ such that for some vertex $t$ of $T$ we have
$v_1\sim_t v_4$ but $v_2\not\sim_t v_4$ and $v_3\not\sim_t v_4$.
Main properties of this subvariety are as follows:

\begin{prop}[\cite{GHP}]\label{properties} 
\begin{enumerate}
\item The codimension of $D(T)$ is equal to the number of the
internal edges (i.e. an edge other than a tail) of $T$.

\item  An inclusion $D(T)\subseteq D(T')$ holds
if and only if $T'$ is obtained by
contracting $T$ along certain internal edges.

\item The subvariety $D(T)$ is canonically isomorphic to $\prod_{t\in
T_0}\overline{{\mathcal M}}_{0,val(t)}$. Here $T_0$ is the set of
vertices of $T$ and $val(t)$ denotes the valency of $t$.

\item Let $D^*(T)=D(T)-\cup D(T')$, where the union is taken over
all the trees that can be contracted to $T$ except $T$ itself.
Then the set of $D^*(T)$ for all non-equivalent labeled trees
gives a stratification of $\MM$.

\item If $T$ and $T'$ are two stable trees with only one internal
edge then the corresponding subvarieties $D(T)$ and $D(T')$ are
codimension one divisors. Let $A_1$ and $A_2$ be the set of tails
attached to the corresponding two vertices of $T$ and similarly
define $A'_1$ and $A_2'$ for $T'$. The divisors $D(T)$ and $D(T')$
intersect if and only if one of the following conditions hold:
$$A_i\subseteq A_j'\quad \mbox{or}\quad A'_i\subseteq A_j$$ for
some $i,j\in\{1,2\}$.
\end{enumerate}
\end{prop}

We consider the following  affine covering of $\MM$. For each
labeled tree $T$ let $U(T)$ be the open subset of $\MM$ given by
$$\lambda_v\ne 0\quad \mbox{for all}\quad v\not\in V(T).$$
Notice that for $v\in V(T)$ we should have $\lambda_v\ne \infty$
on $U(T)$. Because if $\lambda_{v_1v_2v_3v_4}=\infty$ then
$\lambda_{v_2v_1v_3v_4}=0$ but $(v_2v_1v_3v_4)\not\in V(T)$. This
shows that $U(T)$ is affine. The basic properties of these subsets
are collected in the following lemma:

\begin{lem}[\cite{GHP}]
\begin{enumerate}
\item If $T$ is a contraction of $T'$ along some internal edges
then $U(T)\subset U(T')$. Therefore $\M\subset U(T)$ for all
choices of $T$.

\item If the combinatorial tree associated to a marked stable
curve $q\in \MM$ is $T$ then $q\in U(T)$. Therefore we have a
covering of $\MM$ by $U(T)$'s.

\item We have the following relation:
$$U(T)=\MM-\cup D(T')$$
where the union is taken over all the labeled trees $T'$ such that
$T$ can not be contracted to $T'$.
\end{enumerate}
\end{lem}

Finally we need to have an inductive way of constructing a
coordinate system on $U(T)$. In fact we can choose $N$ elements of
$v_1,\dots,v_N\in V_N$ such that the corresponding functions
$\lambda_{v_i}$ form a coordinate system on $U(T)$, i.e. the
coordinate ring of $U(T)$ will be a localization of the polynomial
algebra $\Q[\lambda_{v_1},\dots,\lambda_{v_N}]$. This is done in
\cite{GHP}\S 3.2. For explaining this some preparation is needed.
If $T$ is a stable $(N+3)$-labeled tree and $i$ is a label, then we
define $T\backslash i$ which is a stable $(N+2)$-labeled tree as
follows. If the valency of the vertex of the tail associated to
$i$ is $\geqslant 4$ then we just remove the tail associated to $i$. If
the vertex of the tail associated to $i$ has valency $3$ then
removing the tail makes the tree unstable. Therefore we also
contract the internal edge which has a common vertex with the tail
associated to $i$. Note that the labels of $T\backslash i$ is from
the set $\{1,\dots,N+3\}-\{i\}$.\par

The median of three distinct labels $i,j$ and $k$ in a labeled
tree $T$ is the unique vertex $t$ such that removing $t$ divides
$i,j$ and $k$ into different connected components.
\begin{lem}\label{coordinates}
Let $T$ be an  $(N+3)$-labeled stable tree. Let $i$ be a label such
that the valency of the vertex of the tail associated to $i$,
denoted by $v_i$, is either $\geqslant 4$ or there are at least two
tails attached to $v_i$ (such a label always exists). Inductively
let $\lambda_{v_1},\dots,\lambda_{v_{N-1}}$ be a coordinate system
on $U(T\backslash i)$. Let $(d_1,d_2,d_3)$ be a distinct triple of
labels in $\{1,\dots,N+3\}\backslash\{i\}$ with the following
property: If the valency $v_i$ is $\geqslant 4$ then the median of
$d_1,d_2,d_3$ should be $v_i$. If the valency of $v_i$ is $3$ then
the median of $d_1,d_2,d_3$ should be the unique vertex connected
by an internal edge to $v_i$ and $d_3$ is a tail attached to
$v_i$. Then
$\lambda_{v_1},\dots,\lambda_{v_{N-1}},\lambda_{d_1d_2d_3i}$ is a
coordinate system for $U(T)$.
\end{lem}
{\bf Proof.} Refer to section 3.2 in \cite{GHP}.\qed

\section{The tangential base point}

In this section we will select a special divisor on $\MM$ and
define a particular line in the normal bundle of this divisor
which is symmetric and will play a crucial role in our definition
of the series regularization of $p$-adic MZV's.\par

Let $1\leqslant i\leqslant N$ be an integer, define $T_i$ and $T'_i$ to be
stable $(N+3)$-labeled trees with one internal edge and two vertices
$v_1$ and $v_2$. The tails attached to $v_1$ in $T_i$ have labels
$\{1,\dots,i+1\}$ and in $T_i'$ have labels $\{2,\dots,i+2\}$. The
tails of $v_2$ in $T_i$ have labels $\{i+2,\dots,N+3\}$ and in
$T_i'$ have labels $\{i+3,\dots,N+3,1\}$. We let $D_i=D(T_i)$ and
$D_i'=D(T_i')$ in the notation of the previous section. Observe
that $D_i$'s intersect at one point denoted by $P$ and similarly
$D_i'$ intersect at one point denoted by $Q$. This follows from
the fact that there is a unique tree, which we denote by $T$, with
$N$ internal edges that can be contracted to $T_i$'s, and
similarly there is a unique tree, which we denote by $T'$, with
$N$ internal edges that can be contracted to $T_i'$'s. In fact $T$
is a binary tree with vertices $v_0,\dots,v_{N}$ where $v_i$ and
$v_{i+1}$ are connected by an internal edge and the tails of $v_0$
are $\{1,2\}$, the tail attached to $v_i$ is $i+2$ for $0<i<N$ and
the tails of $v_N$ are $\{N+2,N+3\}$.

\begin{lem}\label{charts}
The collection of $\lambda_{1,N+3,i+2,i+1}$ for $i=1,\dots,N$
gives a coordinate system for $U(T)$. Similarly the collection
$\lambda_{2,1,i+3,i+2}$ for $i=1,\dots, N$ gives a coordinate
system for $U(T')$.
\end{lem}
{\bf Proof.} We only prove the first part, the second part is
similar. The proof is by induction on $N$. For $N=1$,
$U(T)={\mathbb P}^1-\{1,\infty\}$, and $\lambda_{1432}$ sends
$(0,x,1,\infty)$ to $x$, so we have the natural coordinate.
Suppose we have proved the lemma for $N-1$. Then $U(T\backslash
2)$ by induction has coordinates $\lambda_{1,N+3,i+1,i+2}$ for
$i=2,\dots,N$. According to lemma \ref{coordinates} we have to add
$\lambda_{d_1,d_2,1,2}$, such that $(d_1,d_2,1)$ has median $v_1$,
the vertex of the tail associated to $3$ in $T$. This can be
achieved if we let $d_1=3$ and $d_2=N+3$. However since
$\lambda_{3,N+3,1,2}=1-\lambda_{1,N+3,3,2}$ we can use
$\lambda_{1,N+3,3,2}$ as an extra coordinate.\qed

Note that for the point $q=(0,x_1\cdots x_N,x_2\cdots
x_N,\dots,x_N,1,\infty)$ the coordinates
$\lambda_{1,N+3,i+2,i+1}(q)=x_i$. Furthermore
\begin{equation}\label{coor}
z_i:=\lambda_{2,1,i+3,i+2}(q)=\frac{1-x_1\dots x_i}{1-x_1\cdots
x_{i+1}}
\end{equation}
where $x_{N+1}:=0$. Also notice that since $(1,N+3,i+2,i+1)\in
V(T)$ so the equation of $D_i$ inside $U(T)$ is
$\lambda_{1,N+3,i+2,i+1}=0$ or more naively $x_i=0$. Similarly the
equation of $D'_i$ inside $U(T')$ is $\lambda_{2,1,i+3,i+2}=0$ or
$z_i=0$. The divisor $D_N'$ given by $z_N=0$ is the {\it{special
divisor}} that will play an important
 role in defining our regularization.

Let $E_N$ be the Zariski open subset of $\M$ defined by:
$$E_N:={\mathbb A}^N-\left\{x_1\dots x_N=0,\quad \prod_{i\in I}x_i=1,\quad
I\subseteq \{1,\dots,N\}\right\}.$$ For a subset $I$ of
$\{1,\dots,N+3\}$ consisting of non-consecutive numbers we have
the divisor $\prod_{i\in I}x_i=1$ inside $\M$. Its closure inside
$\MM$ is denoted by $D(I)$.

\begin{lem}\label{R}
The line $(t,t,\dots,t)$ inside $\M$ has a limit in $\MM$ when $t$
approaches 1, we denote this point by $R$. This point lies on the
divisor $D_N'$ given by $z_N=0$. Its coordinates using the
coordinate system $z_1,\dots,z_N$ is given by
$(\frac{1}{2},\frac{2}{3},\dots,\frac{N-1}{N},0)$. The point $R$
does not lie on any other component $\MM-E_N$.
\end{lem}
{\bf Proof.} The divisors of $\MM-\M$ are in one to one
correspondence with unordered partitions of $\{1,\dots,N+3\}$ into
two subsets, where each subset has at least two elements. For a
given partition $A\cup B$, the equation of the divisor associated
to it inside $({\mathbb P}^1)^{V_N}$ is given by
$\lambda_{v_1v_2v_3v_4}=0$ for all quadruples such that the sets
$\{v_1,v_4\}$ and $\{v_2,v_3\}$ are separated by $A$ and $B$. Now
since the cross ratio of $(0,t^k,t^l,t^i)$ has the limit
$\frac{l-k}{l-i}$ when $t$ approaches 1 and the cross ratio
$(0,t^k,\infty,t^i)$ has the limit 1, it follows that the limit of
$(0,t^N,\dots,t,1,\infty)$ when $t$ approaches 1 will not lie on
any divisor other than the one obtained by the partition
$\{1,N+3\}\cup\{2,\dots,N+2\}$. In fact since the cross ratio of
$(0,t^i,t^j,\infty)$ approaches 0, the limit point lies on this
divisor. This also shows that $R$ belongs to $U(T')$ using the
above notation. The coordinates of the point
$(0,t^N,\dots,t,1,\infty)$ in terms of $z_i$ are:
$$z_i=\frac{1-t^i}{1-t^{i+1}}\quad (i<N), z_N=1-t^N$$
so when $t$ approaches 1 we get the desired coordinates of the
lemma. If $I$ is a subset of $\{1,\dots,N+3\}$ then the divisor
$D(I)$ in the coordinates $z_i$ is inside the divisor
\begin{equation}\label{zdiv}
1-\prod_{i\in I}\frac{1-z_i\dots z_N}{1-z_{i-1}\dots z_N}=0.
\end{equation}
An easy inspection shows that the above divisor has $z_n=0$ as a
component, if we remove this component then the point $R$ does not
lie on the remaining components. The reason is that substituting
$z_i=\frac{i}{i+1}$ for $i<N$ in \eqref{zdiv} we get:
$$\frac{\prod_{i\in I}(1-\frac{i-1}{N}z_N)-\prod_{i\in I}(1-\frac{i}{N}z_N)}
{\prod_{i\in I}(1-\frac{i-1}{N}z_N)},$$
the numerator is $\frac{|I|}{N}z_N+\cdots$ where the remaining
factors are divisible by $z_N^2$. Now if we divide by $z_N$ and
let $z_N=0$ we get $\frac{|I|}{N}$.\qed

Let ${\mathcal N}^{00}(N)$ be the normal bundle of $D_N'$ minus
the zero section and restricted to $D_N'-\cup D$ where $D$ runs
over all the divisors of $\MM-\M$ other than $D_N'$. According to
the lemma above we have the following embedding:
$$\iota_N:{\mathbb G}_m\hookrightarrow {\mathcal N}^{00}(N)$$
$$t\mapsto\left(\frac{1}{2},\dots,\frac{N-1}{N},Nt\right).$$
\begin{lem}\label{compatible}
The composition
$${\mathbb
G}_m\stackrel{\iota_N}{\To}{{\mathcal
N}^{00}(N)}\stackrel{\pi_N}{\To}{\mathcal N}^{00}(N-1)$$ is
$\iota_{N-1}$.  Here $\pi_N$ is the projection induced from
$\MM\to \overline{\mathcal M}_{0,N+2}$ obtained by neglecting
$(N+2)$-nd marked point.
\end{lem}
{\bf Proof.} Note that $\pi_N$ sends the special divisor of $\MM$
to the special divisor of $\overline{\mathcal M}_{0,N+2}$. Using
the coordinates $z_i$'s around the most exceptional divisors the
equation of $\pi_N$ will become
$(z_1,\dots,z_N)\To (z_1,\dots,z_{N-2},z_{N-1}\cdot z_N)$
(This can be easily derived from the obvious description of the
map in $x_i$ coordinates which is $(x_1,\dots,x_N)\To
(x_1,\dots,x_{N-1})$ and a change of variables). Therefore the
point $\left(\frac{1}{2},\dots, \frac{N-1}{N},Nt\right)$ will map
to $\left(\frac{1}{2},\dots,\frac{N-2}{N-1},(N-1)t\right)$ which
is by definition $\iota_{N-1}(t)$.\qed

\section{Series regularization of $p$-adic multiple zeta values}
Let $1\leqslant p\leqslant m$. We will interpret the MPL:

$$Li_{n_1,\dots,n_m}(x_p,\dots,x_m)=\sum_{0<k_1<\dots<k_m}
\frac{x_p^{k_p}\cdots x_m^{k_m}}{k_1^{n_1}\cdots\cdots k_m^{n_m}}$$
as a Coleman function on $\M$ where $N\geqslant m$. Notice
that the number of variables could be smaller than the depth $m$.
This will be done inductively using the differential equation:

$$dLi_{n_1,\dots,n_m}(x_p,\dots,x_m)=\sum_{i=p}^m
\partial_iLi_{n_1,\dots,n_m}(x_p,\dots,x_m)$$
where $\partial_iLi_{n_1,\dots,n_m}(x_p,\dots,x_m)$ is given by
the following formula
$$
\begin{cases}
Li_{n_1,\dots,n_i-1,\dots,n_m}(x_p,\dots,x_m)d\log x_i
&\text{if}\quad n_i>1,
\\
Li_{n_1,\dots,\widehat{n_i},\dots,n_m}(x_p,\dots,x_{i-1}x_i,\dots,x_m)d\log
(1-x_i)
\\
\quad\quad
-Li_{n_1,\dots,\widehat{n_i},\dots,n_N}(x_p,\dots,x_ix_{i+1},\dots,x_m)d\log
x_i(1-x_i)&\text{if}\quad n_i=1
\end{cases}
$$
where by convention $x_{p-1}=1$ and in the case where $n_m=1$ and
$i=m$ we omit the last line, i.e. formally let $x_{m+1}=0$ and
assume $Li_{n_1,\dots,n_m}(x_p,\dots,x_{m-1},0)=0$. We remark that
the way we have parameterized $\M$ is specially useful to see that
the multiple polylogarithm is a Coleman function.
We denote the corresponding Coleman function by 
$Li^a_{n_1,\dots,n_m}(x_p,\dots,x_m)$
\footnote{We may sometimes omit \lq$a$'.}
in accordance with a branch $a\in\bold Q_p$ of the $p$-adic logarithm.

Note that if instead of the coordinates $x_1,\dots,x_N$ we had
used the permuted coordinates $x_{\tau(1)},\dots,x_{\tau(N)}$ for
$\tau\in{\frak S}_N$,
we get the following diagram.

\begin{equation}\label{diagram}
\begin{CD}
 \mbox{Col}^a(E_N) @<<<  \mbox{Col}^a(\M)@>>>
\mbox{Col}^a({\mathcal N}^{00}) @>\iota_N^{a*}>> \mbox{Col}^a({\mathbb G}_m)\\
@AA\tau^a A                      @AA\tau^a A @AA\tau^a A @AA= A \\
\mbox{Col}^a( E_N^{\tau}) @<<< \mbox{Col}^a(\M^{\tau})@>>>
\mbox{Col}^a(({\mathcal N}^{00})^{\tau}) @>\iota_N^{a*}>>
\mbox{Col}^a({\mathbb G}_m)
\end{CD}
\end{equation}
where $E_N$ is defined by
$$E_N:={\mathbb A}^N-\left\{x_1\dots x_N=0,\quad \prod_{i\in I}x_i=1,\quad
I\subseteq \{1,\dots,N\}\right\}$$ and $\M^{\tau}$ is defined by
$\M^{\tau}=\M\times_{\mathbb A^N,\tau} {\mathbb A}^N$
where $\tau:{\mathbb A}^N\to{\mathbb A}^N$ is a map sending 
$x_i\mapsto x_{\tau(i)}$.
The isomorphism
$id\times\tau:\M\To\M^{\tau}$ induces an isomorphism $\tau:\MM\To
\MM^{\tau}$ which defines $\tau:{\mathcal N}^{00}\To ({\mathcal
N}^{00})^{\tau}$.

 The commutativity of all the squares
follows from the functorial property of analytic continuation and
restriction. The commutativity of the farthest right square
follows from the fact that the line $M=\{(t,t,\dots,t)\}$ is
invariant under the permutation of the variables. We remark that
$\iota_N$ does not have a good reduction modulo $p$ if $p\leqslant N$ and
the definition of $\iota_N$ requires using the approach of
Vologodsky.
\\
 We are now ready to give the
following definition for the series regularized $p$-adic MZV,
which is more accurately an element of $\Q_p[T]$. Its definition a
priori depends on the choice of a branch of $p$-adic logarithm.
i.e. a choice of $a\in \Q_p$ for the value of this logarithm at
$p$.
\begin{defn}
 The series regularized $p$-adic MZV $\zeta_p^S(n_1,\dots,n_m)\in\Q_p[T]$ is
 defined as follows. The differential equations of MPL shows that
 the image of $Li_{n_1,\dots,n_m}(x_1,\dots,x_m)$
 under the following maps
\begin{equation}\label{tangential}
 t^a_{\M}:\mbox{Col}^a(\M)\To \mbox{Col}^a({\mathcal
 N}^{00})\stackrel{\iota_N^{a*}}{\To}\mbox{Col}^a({\mathbb G}_m)
\end{equation}
 is an element of $\Q_p[T]
 \subset{\mathcal O}_{{\mathbb G}_m}[T]=\mbox{Col}^a({\mathbb G}_m)$
 with $T=\log^a t$. The series regularized value is defined to be
 this polynomial.
\end{defn}
This is independent of the choice of $N\geqslant m$ which follows from
the lemma \ref{compatible}. This regularization will be
independent of the choice of the branch $a$. But this will be
proved later in Theorem \ref{branchfree}.

\section{Series shuffle relation}

We now describe the series shuffle relation for multiple
polylogarithms. To do this we need the notion of generalized
shuffles of order $r$ and $s$, denoted by
\begin{align*}
Sh^{\leqslant}(r,s):=\underset{N}{\bigcup}\Bigl\{
\sigma:&\{1,\cdots,r+s\}\to\{1,\cdots,N \}\Bigl| \  \sigma{\text{ is onto}}, \\
&\sigma(1)<\cdots<\sigma(r), \sigma(r+1)<\cdots<\sigma(r+s)
\Bigr\}.
\end{align*}
We recall the definition from \cite{G1}\S 7.1. Let
$$\Z_{++}^m=\{(k_1,\dots,k_m)\in \Z_+^m\quad|\quad
0<k_1<\dots<k_m\}\subset \Z^m.$$ 
There is a natural decomposition:
$$\Z_{++}^r\times\Z_{++}^s=\bigcup_{\sigma\in
Sh^{\leqslant}(r,s)} \Z^{\sigma}_{++},$$ where
$$\Z^{\sigma}_{++}:=\Bigl\{(k_1,\cdots,k_{r+s})\in\Z^{r+s}_{++}\Bigm|
k_i<k_j \text{ if } \sigma(i)<\sigma(j), \
k_i=k_j \text{ if } \sigma(i)=\sigma(j)
\Bigr\}.
$$
For example for $r=s=1$ we have:
$$\{k_1>0\}\times \{k_2>0\}=\{0<k_1<k_2\}\cup
\{0<k_1=k_2\}\cup\{0<k_2<k_1\}.$$ We define the permuted multiple
polylogarithm by:
$$Li_{n_1,\dots,n_{r+s}}^{\sigma}(x_1,\dots,x_{r+s})=
\sum_{(k_1,\dots,k_{r+s})\in\Z_{++}^{\sigma}}\frac{x_1^{k_1}\cdots
x_{r+s}^{k_{r+s}}}{k_1^{n_1}\cdots k_{r+s}^{n_{r+s}}}.$$ Then
formally we have:
\begin{equation}\label{sershuffle}
\begin{split}
Li_{n_1,\dots,n_r}&(x_1,\dots,x_r)Li_{n_{r+1},\dots,n_{r+s}}(x_{r+1},\dots,x_{r+s}) \\
&={\underset{\sigma\in Sh^{\leqslant}(r,s)}{\sum}}
Li^{\sigma}_{n_1,\dots,n_{r+s}}(x_1,\dots,x_{r+s}).
\end{split}
\end{equation}
In fact $Li^{\sigma}_{n_1,\dots,n_{r+s}}(x_1,\dots,x_{r+s})$ is of
the form $Li_{q_1,\dots,q_l}(y_1,\dots,y_l)$ with the same weight
and $y_i$ are either one of the $x_j$'s or product of two of the
$x_j$'s. If we let $N=r+s$, notice that this function can
{\it{not}} be considered in general as a Coleman function on $\M$.
This follows from the fact that the parameterization of $\M$ is
not symmetric with respect to the permutation of the coordinates
$x_1,\dots,x_N$. Recall that we only remove the product of the
consecutive coordinates equalling 1. However all of these
functions can be considered as Coleman functions on $E_N$ defined
in the diagram \ref{diagram}.

Since the shuffle formula given above is formal, its validity can
be extended if we regard the functions as Coleman functions on
$E_N$. The idea of the proof of series shuffle relation is to
analytically continue the formula \eqref{sershuffle} for Coleman
functions on $E_N$ to ${\mathcal N}^{00}$. Then show that if we
restrict it to ${\mathbb G}_m$ via the embedding given before we
get the appropriate relation between the regularized $p$-adic
MZV's. The crucial step is the following proposition which was
inspired from Proposition7.7 of \cite{G1}.
\\
\begin{prop}\label{main}
 Let $1=p_1\leqslant p_2\leqslant \dots\leqslant p_m \lneqq p_{m+1}=N+1$ be integers and
let $y_i=\prod_{j=p_i}^{p_{i+1}-1}x_j$.\footnote{Empty product is
defined to be $1$.} Assume below that $n_k>1$ and $k+l=m$ . Let
$F$ be the following Coleman function on $\M$: for $l>0$, $F$ is:
$$Li^a_{n_1,\dots,n_k,\underbrace{1,\dots,1}_l}(y_1,\dots,y_k,y_{k+1},\dots,y_{k+l})
-Li^a_{n_1,\dots,n_k,\underbrace{1,\dots,1}_l}(y_{k+1},\dots,y_{k+l})$$
and for $l=0$, $F$ is
$$Li^a_{n_1,\dots,n_k}(y_1,\dots,y_{k})-Li^a_{n_1,\dots,n_k}(y_{k}).$$
For a divisor $D$ in $\MM-\M$ let $F^{(D)}$ denote the extension
of $F$ to the normal bundle ${\mathcal N}_D^{00}$. Then
$F^{(D)}=0$ for $D=D_N',D_{N-1}'$ and $D_N$.
\end{prop}
{\bf Proof.} Notice that $D_N'$ intersects $D_{N-1}'$ and
$D_{N-1}'$ intersects $D_N$. In fact we saw in section 2 that all
the divisors $D_i'$ intersect at a single point $Q$. To see that
$D_{N-1}'$ and $D_N$ intersect we can use part (5) of the
proposition
\ref{properties}. \\
We will prove below that $F^{(D_N)}$ is zero and $F^{(D_N')}$ and
$F^{(D_{N-1}')}$ are constant. Now since $D_{N-1}'$ and $D_N$
intersect if we apply Proposition 2.6 of \cite{BF} it follows that
$F^{(D_{N-1}')}$ is zero. A similar argument using the divisors
$D_{N-1}'$ and $D_N'$ implies that $F^{(D_N')}$ vanishes as well.

\begin{lem} The extension of $Li^a_{n_1,\dots,n_m}(y_1,\dots,y_m)$ when
extended to ${\mathcal N}_{D_N}^{00}$ is zero.
\end{lem}
{\bf Proof.} The constant term of MPL at the origin, i.e. the
intersection of $D_i$'s for $i=1,\dots,N$ is zero. This follows
from the fact that in the neighborhood of the origin we have the
power series expansion without constant term. We now calculate the differential of the
MPL and take its residue at $x_N=0$. This will be zero since each
term will be a MPL of weight one smaller than the original MPL and
hence by induction will be zero. This finishes the proof.\qed

Let us now show that $F^{(D_j')}$ is constant for $j=N-
 1,N$.
Recall that the coordinates $z_i$'s for the divisors $D_i'$ are
related to the original coordinates $x_i$ by:
$$x_i=\frac{1-z_i\dots z_N}{1-z_{i-1}\dots z_N}.$$
The residue of $d\log x_i$ at $z_{N}=0$ and $z_{N-1}=0$ is zero
for $i<N$. This implies that the differentials with respect to
those indices $i$ for which $n_i>1$ do not contribute. If $n_i=1$
for $i<k$ then the differential of $F$ with respect to $y_i$ (we
are assuming $y_i\not\equiv 1$) is:
\begin{align*}
&Li_{n_1,\dots,\hat{n_i},\dots,
n_k,\underbrace{1,\dots,1}_l}(y_1,\dots,y_{i-1}y_i,\dots,y_{k+l})d\log
(1-y_i)
\\
\quad \quad
&-Li_{n_1,\dots,\hat{n_i},\dots,n_k,\underbrace{1,\dots,1}_l}(y_1,\dots,y_iy_{i+1},\dots,y_{k+l})d\log
y_i(1-y_i)
\end{align*}
Since $d\log y_i$ has residue zero along $z_N=0$ or $z_{N-1}=0$
an induction on the weight shows that this difference is zero
when the residue is taken. So the only variable that are left are
those $y_{k+1},\dots,y_{k+l}$ that are not identically 1 when
$l>0$ and $y_k$ when $l=0$. The induction implies that these also
do not have any contribution. We provide the details for the case
when $l>0$, the case $l=0$ is similar and even simpler. The
differential with respect to $y_i$ when $i>k$ is given by:
\begin{align*}
&Li_{n_1,\dots,
n_k,\underbrace{1,\dots,1}_{l-1}}(y_1,\dots,y_{i-1}y_i,\dots,y_{k+l})d\log
(1-y_i)
\\
&-Li_{n_1,\dots,
n_k,\underbrace{1,\dots,1}_{l-1}}(y_{k+1},\dots,y_{i-1}y_i,\dots,y_{k+l})d\log
(1-y_i)
\\
&-Li_{n_1,\dots,
n_k,\underbrace{1,\dots,1}_{l-1}}(y_1,\dots,y_{i}y_{i+1},\dots,y_{k+l})d\log
y_i(1-y_i)
\\
&+Li_{n_1,\dots, n_k,
\underbrace{1,\dots,1}_{l-1}}(y_{k+1},\dots,y_{i-1}y_i,\dots,y_{k+l})d\log
y_i(1-y_i)
\end{align*}
Now it is clear that induction on the weights imply that the first
two and the last two terms will cancel each other after taking the
residues. This finishes the proof of the proposition
\ref{main}.\qed

\begin{cor}\label{admissible}
With notation of proposition \ref{main} and $n_k>1$, the analytic
continuation of $Li^a_{n_1,\dots,n_k}(y_1,\dots,y_k)$ to ${\mathcal
N}_{D_N'}^{00}$ is constant and is equal to
$\zeta_p(n_1,\dots,n_k)$.
\end{cor}
{\bf Proof.} By proposition \ref{main} the analytic continuation
is the same as the analytic continuation of
$Li^a_{n_1,\dots,n_k}(y_k)$. If $y_k=x_m$ then the claim follows
from lemma \ref{compatible} and the definition of $p$-adic MZV. If
$y_k=x_i\cdots x_m$,  a similar argument as above using the
differential equation of MPL show that the analytic continuation
of $Li^a_{n_1,\dots,n_k}(y_k)-Li^a_{n_1,\dots,n_k}(x_m)$ to $D_N'$ and
$D_{N-1}'$ is constant and it is zero if it is continued to $D_N$.
This finishes the proof.\qed
\section{Proof of main theorems}

  We would like to deduce the validity of series shuffle relation
  from the equation \eqref{sershuffle}. Recall the notations of the
  diagram 4.1 in section 4. We need to construct a tangential
  morphism $\mbox{Col}^a(E_N)\To \mbox{Col}^a(\mathbb{G}_m)$ that
  extends the morphisms $\mbox{Col}^a(\M)\To
  \mbox{Col}^a(\mathbb{G}_m)$ and is invariant under the action of
  the symmetric group ${\frak S}_N$.\par
  
  This can be achieved by the following argument. Consider
  $$E_N\subset \M\subset \MM.$$
The complement of $E_N$ in $\MM$ is union of certain divisors
$\{D_i\}_{i\in I}$ (not necessarily with normal crossings). Let
$X$ be the smooth variety $\MM\backslash \bigcup_{D_i\neq D_N'}
D_i$. The complement of $E_N$ in $X$ is a Zariski open subset of
$D_N'$ denoted by $D^{\circ}$. By lemma \ref{R}, the point $R$
belongs to this complement and therefore the map
$\iota_N:\mathbb{G}_m\To {\mathcal N}^{00}$ factors through the
normal bundle of $D^{\circ}$ minus its zero section. Therefore the
tangential morphism $$\mbox{Col}^a(E_N)\To \mbox{Col}^a({\mathcal
N}_{D^{\circ}}^{00})$$ can be composed with the restriction map
$\iota_N^{a*}:\mbox{Col}^a({\mathcal
N}_{D^{\circ}}^{00})\To\mbox{Col}^a({\mathbb G}_m)$ to give the
desired tangential morphism:
$$t^a_N:\mbox{Col}^a(E_N)\To \mbox{Col}^a(\mathbb{G}_m).$$

\begin{lem}\label{symmetric}
The map $t_N$ is invariant under the action of the symmetric group
${\frak S}_N$.
\end{lem}

{\bf Proof.} This simply follows from the naturality of the
tangential morphism as was explained for commutativity of the
diagram \eqref{diagram}. Let $\tau\in{\frak S}_N$, we have the
following commutative diagram:
\begin{equation*}
\begin{CD}
 \mbox{Col}^a(E_N) @>{t^a_N}>>  \mbox{Col}^a({\mathcal N}^{00}_{D^{\circ}})\\
@AA\tau A                      @AA\tau A\\
\mbox{Col}^a( E_N^{\tau}) @>(t_N^a)^\tau>> \mbox{Col}^a(({\mathcal
N}_{D^{\circ}}^{00})^{\tau}).
\end{CD}
\end{equation*}
This together with the following commutative diagram of spaces
finishes the proof:
\begin{equation*}
\begin{CD}
 {\mathcal N}_{D^{\circ}}^{00} @<\iota_N<< {\mathbb G}_m\\
@VV\tau V                       @VV= V \\
({\mathcal N}_{D^{\circ}}^{00})^{\tau} @<\iota_N^{\tau}<< {\mathbb
G}_m.
\end{CD}
\end{equation*}
\qed

 Armed with this tangential morphism we now prove the shuffle
relation for regularized $p$-adic MZV's.
\\
\begin{lem}
With notation of the formula \eqref{sershuffle}, the image of
$Li^{\sigma}_{n_1,\dots,n_{r+s}}(x_1,\dots,x_{r+s})$ under the
tangential morphism $t_{r+s}$ is
$\zeta_p^S(\sigma(n_1,\dots,n_{r+s}))$. Here
$\sigma(n_1,\cdots,n_{r+s})=(c_1,\cdots,c_{N})$ where $N$ is the
cardinality of the image of $\sigma$ and
$$c_i=\begin{cases}
n_m+n_l & \text{if  } \sigma^{-1}(i)=\{m,l\}, \\
n_m    & \text{if  } \sigma^{-1}(i)=\{m\}.\\
\end{cases}
$$
\end{lem}
{\bf Proof.} Note that
$$Li_{n_1,\dots,n_{r+s}}^{\sigma}(x_1,\dots,x_{r+s})=Li_{c_1,\dots,c_N}(y_1,\dots,y_{N})$$
where
$$y_i=\begin{cases}
x_mx_l & \text{if  } \sigma^{-1}(i)=\{m,l\}, \\
x_m    & \text{if  } \sigma^{-1}(i)=\{m\}.\\
\end{cases}
$$
Now since the tangential morphism $t^a_{r+s}$ is invariant under the
action of the symmetric group, with a permutation of the
parameters we can assume that we are in the situation of the
proposition \ref{main}. Therefore if $c_N>1$ the result follows
from corollary \ref{admissible}. If for some $k$ we have
$c_{k-1}>1$ and $c_{k}=\dots=c_N=1$ then we necessarily have
$y_k=x_i,\dots,y_{N}=x_{r+s}$ with $i=r+s+k-N$. The proposition
\ref{main} implies that the extension of
$Li_{c_1,\dots,c_{k-1},\underbrace{1,\dots,1}_{N-k+1}}(y_1,\dots,y_N)$
to ${\mathcal N}_{D_{r+s}'}^{00}$ is the same as the extension of
$Li_{c_1,\dots,c_{k-1},\underbrace{1,\dots,1}_{N-k+1}}(x_i,\dots,x_{r+s})$
to ${\mathcal N}_{D_{r+s}'}^{00}$. Another application of the same
proposition implies that this extension is the same as analytic
continuation of
$Li_{c_1,\dots,c_{k-1},\underbrace{1,\dots,1}_{N-k+1}}(x_1,\dots,x_{k-1},x_i,\dots,x_{r+s})$
to ${\mathcal N}_{D_{r+s}'}^{00}$. This is a Coleman function on
$E_{r+s}$ with variables $(x_1,\dots,x_{r+s})$. Using lemma
\ref{symmetric} we can write it in a more standard way as
$Li_{c_1,\dots,c_{k-1},\underbrace{1,\dots,1}_{N-k+1}}(x_1,\dots,x_N)$.
By lemma \ref{compatible} the image of this function under the map
$t^a_{r+s}$ is the same if we consider it as a function on $E_N$ and
apply the map $t^a_N$. By definition therefore the image under $t^a_N$
of this function as a Coleman function on $E_N$ is
$\zeta_p^S(c_1,\dots,c_{N})$.\qed

\begin{thm}
Series shuffle relations for series regularized $p$-adic MZV holds, i.e.
$$
\zeta_p^S(n_1,\cdots,n_r)\cdot\zeta_p^S(n_{r+1},\cdots,n_{r+s})
={\underset{\sigma\in Sh^{\leqslant}(r,s)}{\sum}}
\zeta_p^S(\sigma(n_1,\cdots,n_{m+p}))
$$
holds for all $r,s,n_1,\cdots,n_{r+s}\geqslant 1$.
\end{thm}
{\bf Proof.} Apply the homomorphism $t^a_{r+s}$ to both sides of the
identity \eqref{sershuffle} and use the previous lemma.\qed

It is interesting that this theorem together with the observation
that $\zeta_p^S(1)=-T$ as a polynomial in $\Q_p[T]$, which follows
from the fact that $Li^a_1(z)=-\log^a(1-z)$ implies that the
definition of the series regularization is independent of the branch of
the $p$-adic logarithm.

\begin{thm}\label{branchfree}
The definition of series regularized $p$-adic MZV
$\zeta_p^S(n_1,\cdots,n_m)$ does not depend on the choice of a
branch $a\in\bold Q_p$ of the $p$-adic logarithm.
\end{thm}

{\bf Proof.} This is clear if $n_m>1$, since by corollary
\ref{admissible} we have
$\zeta_p^S(n_1,\dots,n_m)=\zeta_p(n_1,\dots,n_m)$, 
which is independent of the branch (cf. \cite{Fu1}). 
Now assume that $n_k>1$ and
$n_{k+1}=\dots=n_{k+l}=1$  where $k+l=m$. Then the series shuffle
relation implies that
\begin{align*}
\zeta_p^S(1)&\zeta_p^S(n_1,\dots,n_k,\underbrace{1,\dots,1}_{l-1})
=l\zeta_p^S(n_1,\dots,n_k,\underbrace{1,\dots,1}_l)\\
&+\zeta_p^S(n_1,\dots,n_{k-1},1,n_k,\underbrace{1,\dots,1}_{l-1})
+\dots+\zeta_p^S(1,n_1,\dots,n_k,\underbrace{1,\dots,1}_{l-1})\\
&+\zeta_p(n_1+1,\dots,n_k,\underbrace{1,\dots,1}_{l-1})+\dots
+\zeta_p^S(n_1,\dots,n_k,1,\dots,1,2).
\end{align*}
Now an induction on $l$ proves the theorem.
\qed

Let us now explain the integral regularization of $p$-adic MZV's.
The one variable MPL
$$
Li_{n_1,\dots,n_m}(z)=\sum_{0<k_1<\dots<k_m}\frac{z^{k_m}}{k_1^{n_1}\dots k_m^{n_m}}
$$
can be viewed as a Coleman function on $\Mi=E_1$. Its image under
the tangential morphism $t^a_1$ is the integral regularization of
the $p$-adic MZV and we use the notation
$\zeta_p^I(k_1,\cdots,k_m)$. This is an element of $\Q_p[T]$ where
$T=\log^a (1-z)$. By the $p$-adic iterated integral expression of
$p$-adic MPL, the first author in \cite{Fu1} deduced an integral
shuffle product formula that we now explain. For $\bold
k=(k_1,\dots,k_m)$ and ${\bold k'}=(k'_1,\dots,k'_{m'})$ with
$k_i, k'_j\geqslant 1$ the following formula holds for $p$-adic
MPL's:
\begin{equation}\label{lishuffle}
Li_{\bold k}(z)Li_{\bold k'}(z)=\underset{\tau\in Sh(N,N')}{\sum}
Li_{{\bold a}_{\tau(W_{\bold k},W_{\bold k'})}}(z).
\end{equation}
Here $N=k_1+\cdots+k_m$, $N'=k'_1+\cdots+k'_{m'}$ and
\begin{align*}
Sh(N,N'):=\Bigl\{\tau:\{1,\dots,&N+N'\}\to\{1,\dots,N+N'\}\Bigm|
\tau\text{ is bijective}, \\
&\tau(1)<\cdots<\tau(N), \tau(N+1)<\cdots\tau(N+N')\Bigr\}
\end{align*}
For $W=X_1\cdots X_k$, $W'=X_{k+1}\cdots X_{k+l}$ with
$X_i\in\{A,B\}$ and $\tau\in Sh(k,l)$, the symbol $\tau(W,W')$
stands for $Z_1\cdots Z_{k+l}$ with $Z_i=X_{\tau^{-1}(i)}$. For
$\bold a=(a_1,\cdots,a_l)$ with $l,a_1,\dots,a_l\geqslant 1$ the
symbol $W_{\bold a}$ means a word $A^{k_l-1}BA^{k_{l-1}-1}B\cdots
A^{k_1-1}B$ and for such $W$ we denote its corresponding index by
${\bold a}_W$.\par Each term in \eqref{lishuffle} lies in
$\mbox{Col}^a(\Mi)$. Applying the morphism $t_{\Mi}^a$ to identity
\eqref{lishuffle} gives the following:

\begin{prop}
The integral series shuffle relation for integral regularized
$p$-adic MZV's holds, i.e.
$$
\zeta_p^I(\bold k)\zeta_p^I(\bold k')=\underset{\tau\in
Sh(N,N')}{\sum}\zeta_p^I({\bold a}_{\tau(W_{\bold k},W_{\bold
k'})})
$$
holds for $\bold k=(k_1,\dots,k_m)$ and ${\bold
k'}=(k'_1,\dots,k'_{m'})$.
\end{prop}

Note that $\zeta_p^I(1)=-T$ and therefore the integral shuffle
relation implies that
\begin{equation}\label{zeta1}
\zeta_p^I(\underbrace{1,\dots,1}_n)=\frac{(-T)^n}{n!}
\end{equation}
The proof of the regularization relation is a $p$-adic analogue of
the proof given in section 7 of \cite{G1}.


\begin{thm}
The regularization relation holds. Namely for
$n_1,\dots,n_m\geqslant 1$
\begin{equation}\label{reg}
\zeta_p^S(n_1,\dots,n_m)={\mathbb L}_p(\zeta_p^I(n_1,\dots,n_m))
\end{equation}
where ${\mathbb L}_p:\Q_p[T]\To \Q_p[T]$ is a linear map that is
defined by:
$$\sum_{n=1}^{\infty}{\mathbb
L}_p(T^n)\frac{u^n}{n!}=\exp\left(-\sum_{n=1}^{\infty}\frac{\zeta_p^I(n)}{n}u^n\right).$$
\end{thm}
{\bf Proof.} The validity of the equation \eqref{reg} is clear if
$n_m>1$.  The following special case for the case when
$(n_1,\dots,n_m)=(1,\dots,1)$ can be proved exactly as in lemma
7.9 of \cite{G1}
\begin{equation}\label{zeta1}
\sum_{n=1}^{\infty}\zeta_p^S(\underbrace{1,\dots,1}_n)u^n=\exp\left(-\sum_{n=1}^{\infty}\frac{\zeta_p^I(n)}{n}(-u)^n\right).
\end{equation}

 Assume that $n_k>1$ and $n_{k+1}=\dots=n_m=1$. We
prove the regularization formula by induction on $m-k$. Note that
\begin{align}\label{1}
\begin{split}
Li_{n_1,\dots,n_k}(&x_1,\dots,x_k)Li_{\underbrace{1,\dots,1}_{m-k}}(x_{k+1},\dots,x_m)\\
&=Li_{n_1,\dots,n_m,\underbrace{1,\dots,1}_{m-k}}(x_1,\dots,x_m)+\mbox{other
terms,}
\end{split}
\end{align}
\begin{align}\label{2}
\begin{split}
Li_{n_1,\dots,n_k}(&\underbrace{1,\dots,1}_{k-1},x)Li_{\underbrace{1,\dots,1}_{m-k}}(\underbrace{1,\dots,1}_{m-k-1},y)\\
&=Li_{n_1,\dots,n_m,\underbrace{1,\dots,1}_{m-k}}(\underbrace{1,\dots,1}_{k-1},x,\underbrace{1.\dots,1}_{m-k-1},y)+\mbox{other
terms.}
\end{split}
\end{align}
We apply the tangential map $t^a_N:\mbox{Col}^a(E_m)\To
\mbox{Col}^a({\mathbb G}_m)$ to the equation \eqref{1} and the
tangential map $t^a_2:\mbox{Col}^a({\mathcal M}_{0,5})\To
\mbox{Col}^a({\mathbb G}_m)$ to the equation \eqref{2}. Now if we
use proposition \ref{main} it follows that the first equation
gives the series regularized $p$-adic MZV's and the second will
give the integral regularized $p$-adic MZV's. We therefore have:
\begin{equation}\label{3}\zeta_p^S(n_1,\dots,n_k)\zeta_p^S({\underbrace{1,\dots,1}_{m-k}})=\zeta_p^S({n_1,\dots,n_m,\underbrace{1,\dots,1}_{m-k}})+\mbox{other
terms,}\end{equation}
\begin{equation}\label{4}
\zeta_p^I(n_1,\dots,n_k)\zeta_p^I({\underbrace{1,\dots,1}_{m-k}})=\zeta_p^I({n_1,\dots,n_m,\underbrace{1,\dots,1}_{m-k}})+\mbox{other
terms.}
\end{equation}
 The left hand side of the equation \eqref{4} after applying the map ${\mathbb L}_p$ coincides with the left hand side of
 \eqref{3}. This follows from the equation \eqref{zeta1}. Also note
 that the terms which are not written in the equations \eqref{4}
 and \eqref{4} have less than $m-k$ one at the end so by induction
 after applying ${\mathbb L}_p$ will match. This finishes the
 proof.\qed

\section{Deligne's problem on double shuffle
relations}\label{deligne}

In \cite{D1}, Deligne proposed the following definition for
$p$-adic MZV's. Let $X={\mathbb P}^1\backslash \{0,1,\infty\}$ and
$\pi^{\text{DR}}(X,\overrightarrow{01})$ denote the de Rham
fundamental group of $X$ with the tangential base point
$\overrightarrow{01}$ at $0$. This can be identified as the group
like elements with constant term 1 of the noncommutative power
series Hopf algebra $\Q\langle\langle A,B\rangle\rangle$, where
$A$ corresponds to the loop around $0$ and $B$ to the loop around
$1$. The coproduct is defined by $\Delta A=A\otimes 1+1\otimes A$
and similarly for $B$. Since $X$ and the base point have a good
reduction modulo $p$, we have an action of the Frobenius
endomorphism $\phi$ on this fundamental group tensored with $\Q_p$
which can be extended to $\Q_p\langle\langle A,B\rangle\rangle$.
It can be shown that under this endomorphism\footnote{We are using
the inverse of the usual Frobenius as opposed to \cite{D1}}:
$$A\mapsto \frac{A}{p}$$
$$B\mapsto (\Phi_{De}^p)^{-1}(\frac{B}{p})\Phi_{De}^p$$
for a certain group like element $\Phi_{De}^p$ of
$\Q_p\langle\langle A,B\rangle\rangle$ with constant term 1 and
the coefficients of $B^n$ equal zero. Deligne defines
$(-1)^m\zeta_p^{De}(n_1,\dots,n_m)$ to be the coefficient of
$A^{n_m-1}B\cdots A^{n_m-1}B$ in $\Phi_{De}^p$.

It was asked in \cite{D1} and \cite{DG} to prove the validity of
generalized double shuffle relation for these $p$-adic MZV's. This
is achieved using the results of this paper and \cite{Fu2}. In
fact in the language of Racinet in \cite{R} we need to show that
$\Phi_{De}^p(A,-B)\in \underline{DMR}_0(\Q_p)$. We recall his
machinery briefly. The group scheme $\underline{DMR}_0$ has $k$
(:a field of characteristic $0$)-valued points which is a subset of
power series $k\langle\langle A, B\rangle\rangle$ of those power
series $g=\sum c_W W$, where $W$ runs over all words in $A$ and
$B$, such that:
\begin{enumerate}

\item The constant term $c_{\emptyset}=1$ and $c_A=c_B=0$.

\item The series $g$ is group like with respect to the coproduct defined
above. i.e. $\Delta g=g\otimes g$. This is a concise way
of saying that the coefficients of $g$ satisfies the integral
shuffle relation.

\item Let $\pi_y:k\langle\langle A, B\rangle\rangle\To
k\langle\langle y_1,y_2,\dots\rangle\rangle$ be defined as a
linear map that sends all the words ending $A$ to zero and the
word $A^{n_1-1}B\cdots A^{n_m-1}B$ to $y_{n_1}\cdots y_{n_m}$.
Define the coproduct $\Delta_*$ on $k\langle\langle
y_1,y_2,\dots\rangle\rangle$ by
$$\Delta_* y_n=\sum_{i=0}^n y_i\otimes y_{n-i}\quad\quad y_0:=1.$$
Also define
$$g_*=\exp\left(-\sum_{n=1}^{\infty}\frac{(-1)^n}{n}
c_{A^{n-1}B}y_1^n\right)\pi_y(g).$$ The last condition is that
$\Delta_* g_*=g_*\otimes g_*$. This is a concise way of
saying that the coefficients of $g_*$ satisfy the series
shuffle relation.
\end{enumerate}

In \cite{Fu1} a fundamental solution, denoted by $G_0(z)(A,B)$,
for the $p$-adic KZ equation:
$$dG(z)=\left(A\frac{dz}{z}+B\frac{dz}{z-1}\right)G(z)\quad\quad z\in {\mathbb P}^1(\C_p)\backslash \{0,1,\infty\}$$
was constructed. Its
coefficients are Coleman functions on ${\mathbb
P}^1\backslash\{0,1,\infty\}$. If we analytically continue this
function to the tangent vector $1$ at $z=1$, we get a power series
$\Phi^p_{KZ}(A,B)\in \Q_p\langle\langle A, B\rangle \rangle$,
called the $p$-adic Drinfel'd associator,
whose coefficient for $A_{n_m-1}B\cdots A^{n_1-1}B$ is the integral
regularization $(-1)^m\zeta_p^I(n_1,\dots,n_m)$ evaluated at
$T=0$. The main result of this paper says that in the language of
Racinet
$$\Phi_{KZ}^p(A,-B)\in \underline{DMR}_0(\Q_p).$$
Now we can prove the following:
\begin{thm}
Deligne's $p$-adic MZV satisfies generalized double shuffle
relation.
\end{thm}
{\bf Proof.} It is shown in theorem 2.7 of \cite{Fu2} that
$$\Phi_{KZ}^p(A,-B)=\Phi_{De}^p(A,-B)\cdot
\Phi_{KZ}^p\left(\frac{A}{p},\Phi_{De}^p(A,-B)^{-1}\frac{B}{p}\Phi_{De}^p(A,-B)\right).$$
This can be rewritten using the product $\circledast$ of 
$\underline{DMR}_0$:
$$\Phi_{KZ}^p(A,-B)=\Phi_{KZ}(\frac{A}{p},-\frac{B}{p})\circledast
\Phi_{De}^p(A,-B).$$ The set $\underline{DMR}_0(\Q_p)$ and
also $\bold Q_p\langle\langle A,B\rangle\rangle$ form a group
under this product (ref. \cite{R}) and the two elements
$\Phi_{KZ}^p(A,-B)$ and $\Phi_{KZ}^p(\frac{A}{p},-\frac{B}{p})$
belong to this group, hence $\Phi_{De}^p(A,-B)\in
\underline{DMR}_0(\Q_p)$. Let
$$\tilde{\Phi}_{De}^p(A,B):=\exp(BT)\Phi_{De}^p(A,B).$$ The
coefficient of $A^{n_m-1}B\cdots A^{n_1-1}B$ in
$\tilde{\Phi}_{De}^p$ is denoted by $
(-1)^m\zeta_p^{De,I}(n_1,\dots,n_m)$. This is
$(-1)^m\zeta_p^{De}(n_1,\dots,n_m)$ if $n_m>1$ and if $n_m=1$ this
is a polynomial in terms of $T$ for which if we let $T=0$ we get
$(-1)^m\zeta_p^{De}(n_1,\dots,n_m)$. The series regularization is
obtained by applying ${\mathbb L}_p$, i.e. the coefficient of
$A^{n_m-1}B\cdots A^{n_1-1}B$ in ${\mathbb
L}_p(\tilde{\Phi}_{De}^p)$ is defined to be
$(-1)^m\zeta_p^{De,S}(n_1,\dots,n_m)$.  The fact that
$\Phi_{De}^p\in \underline{DMR}_0(\Q_p)$ implies that
$\zeta_p^{De,I}(n_1,\dots,n_m)$ will satisfy the integral shuffle
relations and $\zeta_p^{De,S}(n_1,\dots,n_m)$ will satisfy the
series shuffle relations. Note that the relation between the two
regularization is automatically holds by the way we have defined
the second regularization.\qed


\end{document}